\documentclass[a4paper,11pt]{amsart}
\usepackage{amssymb}
\newtheorem{theorem}{Theorem}
\newtheorem{thm}[theorem]{Theorem}
\newtheorem{lemma}[theorem]{Lemma}
\newtheorem{prop}[theorem]{Proposition}

\theoremstyle{definition}

\theoremstyle{remark}
\newtheorem{remark}[theorem]{Remark}
\newcommand{\DD}{{\mathbb D}}
\newcommand{\OO}{{\mathcal O}}

\newcommand{\MM}{{\mathcal M}}

\newcommand{\CC}{{\mathbb C}}
\newcommand{\cC}{{\mathcal C}}

\newcommand{\GG}{{\mathbb G}}

\newcommand{\dif}[2]{\frac{\partial #1}{\partial #2}}
 
 \DeclareMathOperator{\diag}{diag}

\DeclareMathOperator{\sgn}{sgn}

\renewcommand{\phi}{\varphi}

\begin{document}
\title
{Geometry of the symmetrized polydisc}

\author{Armen Edigarian}

\address{Institute of Mathematics, Jagiellonian University,
Reymonta 4, 30-059 Kra\-k\'ow, Poland}
\email{Armen.Edigarian@im.uj.edu.pl}
\thanks{The authors were supported in part by the KBN grant
No. 5 P03A 033 21}

\author{W\l odzimierz Zwonek}

\email{zwonek@im.uj.edu.pl}





\begin{abstract} We describe all proper holomorphic mappings of the
symmetrized polydisc and study its geometric properties. We also
apply the obtained results to the study of the spectral unit ball
in $\MM_n(\CC^n)$.
\end{abstract}

\maketitle


\section{Introduction}
Let $\DD$ be the unit disc in the complex plane $\CC$. Let
$\pi_n=(\pi_{n,1},\dots,\pi_{n,n}):\CC^n\to\CC^n$, $n\ge1$, be
defined as follows
\begin{equation*}
\pi_{n,k}(\lambda_1,\dots,\lambda_n)= \sum_{1\le j_1<\dots<j_k\le
n}\lambda_{j_1}\dots\lambda_{j_k},\quad k=1,\dots,n.
\end{equation*}

Observe that $\pi_n$ is a proper holomorphic mapping
with the multiplicity equal to $n!$ (see e.g. \cite{Rud}).
Moreover, $\pi_n^{-1}\big(\pi_n(\DD^n)\big)=\DD^n$.
Hence, $\pi_n|_{\DD^n}:\DD^n\to\pi_n(\DD^n)$ is a proper
holomorphic mapping. Put $\GG_n=\pi_n(\DD^n)$ and
$\delta_n=\pi_n\big((\partial\DD)^n\big)$. The domain $\GG_n$ is
called the {\it symmetrized $n$-disc}.

Below we present a number of
results on the geometry of the symmetrized polydisc $\GG_n$.

One of the main results of the paper is to give the following
characterization of proper holomorphic self-mappings of the
symmetrized polydisc.
\begin{theorem}\label{thm:1} Let $f:\GG_n\to\GG_n$ be a holomorphic mapping. Then
$f$ is proper if and only if there exists a finite Blaschke
product $B$ such that
\begin{equation*}
f\big(\pi_n(\lambda_1,\dots,\lambda_n)\big)=
\pi_n\big(B(\lambda_1),\dots,B(\lambda_n)\big),\quad
\lambda_1,\dots,\lambda_n\in\DD.
\end{equation*}
In particular, $f$ is an automorphism if and only if
\begin{equation*}
f\big(\pi_n(\lambda_1,\dots,\lambda_n)\big)=
\pi_n(h(\lambda_1),\dots,h(\lambda_n)\big), \quad
\lambda_1,\dots,\lambda_n\in\DD,
\end{equation*}
where $h$ is an automorphism of the unit disc $\DD$.
\end{theorem}
Note that, if $\psi:\DD\to\DD$ is a holomorphic function, then the
mapping $f_\psi:\GG_n\to\GG_n$, defined as
\begin{equation*}
f_\psi\big(\pi_n(\lambda_1,\dots,\lambda_n)\big)=
\pi_n\big(\psi(\lambda_1),\dots,\psi(\lambda_n)\big),
\quad\lambda_1,\dots,\lambda_n\in\DD,
\end{equation*}
is a well-defined holomorphic mapping. Moreover, $f_\psi$ is
proper (resp.~an automorphism) if and only if $\psi$ is proper
(resp.~an automorphism).

We get Theorem~\ref{thm:1} as a corollary of the following
\begin{theorem}\label{thm:2} Let $f:\DD^n\to\GG_n$ be a holomorphic mapping.
Then $f$ is proper if and only if there exist finite Blaschke
products $B_1$, \dots, $B_n$ such that
\begin{equation*}
f(\lambda_1,\dots,\lambda_n)=\pi_n\big(
B_1(\lambda_1),\dots,B_n(\lambda_n)\big), \quad
\lambda_1,\dots,\lambda_n\in\DD.
\end{equation*}
\end{theorem}

The symmetrized polydisc has been recently studied by many
authors, especially in two-dimensional case (see
e.g.~\cite{Agl-You1, Agl-You2, Agl-You3}, \cite{Cos1},
\cite{Pfl-Zwo}). Description of automorphisms in $\GG_2$ was given
in \cite{Jar-Pfl2}, description of proper mappings in $\GG_2$ was
given in \cite{Edi}. We follow the ideas of the latter paper.


\section{Proofs}

For the set $K=\{j_1,\ldots,j_k\}\subset I_n:=\{1,\ldots,n\}$,
$1\leq j_1<\ldots<j_k\leq n$ and
$\lambda=(\lambda_1,\ldots,\lambda_n)\in\CC^n$ define
$\lambda_K:=(\lambda_{j_1},\ldots,\lambda_{j_k})$.

Let $\lambda_0\in\CC$, $0\le k\le n$. We define $D_k^n(\lambda_0)$
to be the domain in $\CC^n$ of all
$\lambda=(\lambda_1,\ldots,\lambda_n)\in\CC^n$ such that there is
a (uniquely defined) set $K=K(\lambda_0)\subset I_n$ with $\#K=k$
and $|\lambda_j-\lambda_0|<|\lambda_k-\lambda_0|$, $j\in K$, $k\in
I_n\setminus K$.

Note that $D_n^n(\lambda_0)=D_0^n(\lambda_0)=\CC^n$ for any
$\lambda_0\in\CC$.

\begin{remark}\label{rem:3} Let us fix $\lambda_0\in\CC$, $0\leq k\leq n$. It is
immediate to see that the mappings
\begin{multline}
\rho_{k,n}:D_k^n(\lambda_0)\owns\lambda\mapsto\pi_k(\lambda_K)\in\CC^k,\\
\tilde\rho_{k,n}:D_k^n(\lambda_0)\owns\lambda\mapsto\pi_{n-k}(\lambda_{I_n\setminus
K})\in\CC^{n-k}
\end{multline}
are holomorphic. Moreover, for any mapping
$\phi\in\OO(D,\pi_n(D_k^n(\lambda_0)))$, where $D$ is a domain in
$\CC^m$, the mappings $\rho_{k,n}\circ\pi_n^{-1}\circ\phi$ and
$\tilde\rho_{k,n}\circ\pi_n^{-1}\circ\phi$ are holomorphic, too.
\end{remark}

\begin{lemma}\label{lem:4} Let $\phi\in\OO(D,\CC^n)$, where $D$ is a domain
in $\CC^m$. Assume that $\phi(D)\subset\delta_n$. Then $\phi$ is
constant.
\end{lemma}
\begin{proof} We use induction on $n$. The case $n=1$ is trivial.
So assume that $n\geq 2$ and that lemma is valid for dimensions
$1,2,\ldots,n-1$. Let $f=(f_1,\ldots,f_n)$ denote the multi-valued
mapping $\pi_n^{-1}\circ\phi$. Let $k$ denote the maximal number
of elements of a set $K\subset I_n$ such that for some
$\mu_0\in\DD$ all the coordinates $f_j(\mu_0)$, $j\in K$, are
equal. Without loss of generality we may assume that
$f_1(\mu_0)=\ldots=f_k(\mu_0)$ and $f_l(\mu_0)\neq f_1(\mu_0)$,
$l\geq k+1$. Shrinking $D$, if necessary, we may assume that
$\phi(D)\subset\pi_{n}(D_k^n(f_1(\mu_0)))$, so in view of
Remark~\ref{rem:3} (with $\lambda_0=f_1(\mu_0)$) we get that
$\rho_{k,n}\circ f$ is holomorphic. In particular,
$h:=(\rho_{k,n}\circ f)_1$ is holomorphic on $D$ and $|h|$ attains
its maximum (equal to $k$) at $0$. Therefore, losing no
generality, we may assume that $f_1,\ldots,f_k$ are constant on
$D$. In the case $k=n$ this finishes the proof. So assume that
$1\leq k<n$. Note that the function
$\tilde\phi:=\tilde\rho_{k,n}\circ f$ is holomorphic on $D$ and
$\tilde\phi(D)\subset\delta_{n-k}$. Then the inductive assumption
implies that $\tilde\phi$ is constant, which easily implies that
$f_{k+1},\ldots,f_n$ are constant, which finishes the proof.
\end{proof}

Let us define one more mapping. For $0\leq k\leq n$ and for
$w=\pi_k(\lambda)$, $z=\pi_{n-k}(\mu)$ define
\begin{equation}
p_{k,n}:\CC^k\times\CC^{n-k}\owns(w,z)\mapsto\pi_n(\lambda,\mu)\in\CC^n.
\end{equation}
Note that $p_{k,n}$ is a holomorphic mapping,
$p_{k,n}(\GG_k\times\GG_{n-k})=\GG_n$, $p_{k,n}(\bar
\GG_k\times\bar\GG_{n-k})=\bar\GG_n$.

\begin{lemma}\label{lem:5}
Let $\phi\in\OO(D,\CC^n)$, where $D$ is a domain in $\CC^m$, be such that
$\phi(D)\subset\partial\GG_n$. Then there are a $k$ with $1\leq
k\leq n$, $\theta\in\delta_k$, and $\psi\in\OO(D,\GG_{n-k})$ such
that $\phi=p_{k,n}\circ(\theta,\psi)$.
\end{lemma}
\begin{proof} The case $n=1$ is trivial. So assume that $n\geq 2$.
Define as in Lemma~\ref{lem:4}
$(f_1,\ldots,f_n):=\pi_n^{-1}\circ\phi$. Put $N(\lambda):=\#\{j\in
I_n:|f_j(\lambda)|=1\}$, $\lambda\in D$. Define
$k:=\max\{N(\lambda):\lambda\in D\}$. Let $k=N(\mu_0)$, $\mu_0\in
D$. Obviously, $1\leq k\leq n$.

Denote $u(\lambda):=\max\sb{1\leq j_1<\ldots<j_k\leq
n}\{|f_{j_1}(\lambda)|\cdot\ldots \cdot|f_{j_k}(\lambda)|\}$,
$\lambda\in D$. Then $u$ is plurisubharmonic in $D$, $u\leq 1$ on
$D$ and $u(\mu_0)=1$. Then $u\equiv 1$ on $D$. Therefore, for any
$\lambda\in D$ there is a set $K\subset I_n$ with $k$ elements
such that $|f_j(\lambda)|=1$, $j\in K$ and $|f_j(\lambda)|<1$,
$j\in I_n\setminus K$. Then Lemma~\ref{lem:4} finishes the proof
in the case $k=n$ (because $\phi(D)\subset\delta_n$). So assume
that $k<n$. Applying Remark~\ref{rem:3} (for $\lambda_0=0$) we see
that $\phi=p_{k,n}\circ(\eta,\psi)$, where $\eta=\rho_{k,n}\circ
f$, $\psi=\tilde\rho_{k,n}\circ f$ are holomorphic on $D$. It
easily follows from the definition that $\psi(D)\subset\GG_{n-k}$.
Moreover, $\eta(D)\subset\delta_k$, so, in view of
Lemma~\ref{lem:4}, $\eta$ is constant.
\end{proof}

Let $\lambda=(\lambda_1,\dots,\lambda_n)\in\CC^n$. Note that
$\lambda_1,\dots,\lambda_n$ are roots of the polynomial equation
\begin{multline}
z^n-\big(\pi_n(\lambda)\big)_1z^{n-1}+ \big(\pi_n(\lambda)\big)_2
z^{n-2}+\dots+(-1)^{n-1}\big(\pi_n(\lambda)\big)_{n-1}z\\
+(-1)^{n}\big(\pi_n(\lambda)\big)_n=0,\quad z\in\CC.
\end{multline}

Therefore, Lemma~\ref{lem:5} easily implies the following
\begin{lemma}\label{lem:6} Let $\phi\in\OO(D,\CC^n)$, where $D$ is a domain in
$\CC^m$, be such that $\phi(D)\subset\partial\GG_n$. Then there is
a constant $C\in\partial\DD$ such that
\begin{multline}
C^n-\phi_1(\lambda)C^{n-1}+\phi_2(\lambda)C^{n-2}+\dots\\
+(-1)^{n-1}\phi_{n-1}(\lambda)C+ (-1)^{n}\phi_n(\lambda)=0,\quad
\lambda\in D.
\end{multline}
\end{lemma}

The proof of Theorem~\ref{thm:2} is based on the following result.
\begin{prop}\label{prop:7} Let $f:\DD^n\to\GG_n$ be a proper holomorphic
mapping. Then there exists a bounded holomorphic function $B$ on
$\DD$ such that
\begin{multline}
B^n(z_n)-B^{n-1}(z_n)f_1(z)+\dots+(-1)^{n-1}B(z_n)f_{n-1}(z)+(-1)^nf_n(z)
=0\\
z=(z_1,\dots,z_n)\in \DD^n.
\end{multline}
Moreover, $B$ is non-constant.
\end{prop}

\begin{proof} We use methods similar to the ones of Remmert-Stein
(see e.g. \cite{N}). Take a sequence $\DD\ni
z_n^{\nu}\to\partial\DD$. Then there exists a subsequence $\nu_k$
such that $f(\cdot,z_n^{\nu_k})\to\phi$ locally uniformly on
$\DD^{n-1}$. Note that $\phi:\DD^{n-1}\to\partial\GG_n$ is a
holomorphic mapping. Hence, in view of Lemma~\ref{lem:6} there
exists a constant $C\in\partial\DD$ such that
\begin{equation}\label{eq:1}
C^n-C^{n-1}\phi_1+\dots+(-1)^{n-1}C\phi_{n-1}+(-1)^n\phi_n\equiv0\quad\text{
on } \DD^{n-1}.
\end{equation}

Note that
\begin{multline}\label{eq:2}
-C^{n-1}\frac{\partial \phi_1}{\partial
z_j}+\dots+(-1)^{n-1}C\frac{\partial \phi_{n-1}}{\partial
z_j}+(-1)^n\frac{\partial \phi_n}{\partial z_j}\equiv0 \quad\text{
on } \DD^{n-1}\\ j=1,2,\dots,n-1.
\end{multline}

Since $f$ is proper, we have
\begin{equation*}
\det\left[\dif{f_i}{z_j}\right]_{i,j=1,\dots,n}\not\equiv0\quad\text{
on }\DD^n.
\end{equation*}
So, there exists a $k\in\{1,\dots,n\}$ such that
\begin{equation*}
\det\left[\dif{f_i}{z_j}\right]_{i=1,\dots,n, j=1,\dots,n-1, i\ne
k}\not\equiv0\quad\text{ on }\DD^n.
\end{equation*}
We want to solve the equations (w.r.t. $C$)
\begin{equation}
\begin{cases}
\sum_{m=1,m\ne k}^n (-1)^m C^{n-m}\dif{\phi_m}{z_1}=(-1)^{k+1}C^{n-k}\dif{\phi_k}{z_1}\\
\dots \dots\\
\sum_{m=1,m\ne k}^n (-1)^m
C^{n-m}\dif{\phi_m}{z_{n-1}}=(-1)^{k+1}C^{n-k}\dif{\phi_k}{z_{n-1}}
\end{cases}.
\end{equation}

Let
\begin{equation*}
\Phi_m(\phi_1,\dots,\phi_n):=
\det\left[\dif{\phi_i}{z_j}\right]_{i=1,\dots,n, i\ne m,
j=1,\dots,n-1}\text{ on $\DD^{n-1}$}.
\end{equation*}
Then $(-1)^m C^{n-m}\Phi_k(\phi_1,\dots,\phi_n)=(-1)^k
C^{n-k}\Phi_m(\phi_1,\dots,\phi_n)$.

If $k>1$ then take $m=k-1$. Hence,
\begin{equation}\label{eq:3}
\Phi_{k-1}(\phi_1,\dots,\phi_n)=-C\Phi_k(\phi_1,\dots,\phi_n)\quad
\text{ on }\DD^{n-1}.
\end{equation}
We have
\begin{multline}\label{eq:17}
\Phi_{k-1}(\phi_1,\dots,\phi_n)
\dif{\Phi_{k}(\phi_1,\dots,\phi_n)}{z_j}=\\
\Phi_{k}(\phi_1,\dots,\phi_n)\dif{\Phi_{k-1}(\phi_1,\dots,\phi_n)}{z_j}
\quad\text{ on }\DD^{n-1}
\end{multline}
for any $j=1,\dots,n-1$. Note that the equations in \eqref{eq:17}
hold for any choice of possible sequences $z_n^\nu\to\partial\DD$.
Therefore, the maximum principle for holomorphic functions implies
that similar equations hold for $f_1,\dots,f_n$. So,
\begin{equation}
\Phi_{k-1}(f_1,\dots,f_n) \dif{\Phi_{k}(f_1,\dots,f_n)}{z_j}=
\Phi_{k}(f_1,\dots,f_n)\dif{\Phi_{k-1}(f_1,\dots,f_n)}{z_j} \text{
on $\DD^n$}
\end{equation}
for any $j=1,\dots,n-1$. Put $A=\{\Phi_k(f_1,\dots,f_n)=0\}$ and
\begin{equation}
B=-\frac{\Phi_{k-1}(f_1,\dots,f_n)}{\Phi_{k}(f_1,\dots,f_n)}.
\end{equation}
Note that $A$ is a proper analytic subset of $\DD^n$ and that $B$
is a holomorphic function on $\DD^n\setminus A$. Moreover,
$\dif{B}{z_j}=0$ on $\DD^n\setminus A$ for any $j=1,\dots,n-1$.
Hence, $B$ depends only on $z_n$.

Moreover, (use \eqref{eq:1}, \eqref{eq:3})
\begin{equation*}
\Phi_{k-1}^n+\phi_1\Phi_{k-1}^{n-1}\Phi_{k}+\dots+\phi_{n-1}\Phi_{k-1}\Phi_{k}^{n-1}+
\phi_n\Phi_{k}^n\equiv0\quad\text{ on }\DD^{n-1}.
\end{equation*}
Hence, similar result holds for $f_1,\dots,f_n$. From this we get
\begin{equation}
B^n(z_n)-B^{n-1}(z_n)f_1+\dots+(-1)^{n-1}B(z_n)f_{n-1}+(-1)^nf_n\equiv0\quad\text{
on } \DD^n\setminus A.
\end{equation}
Note that $B$ is a bounded function on $\DD^n\setminus A$, so it
extends holomorphically to $\DD^n$.

If $k=1$ then take $m=2$. Later we prove in a similar way.
\end{proof}

\begin{proof}[Proof of Theorem~\ref{thm:1}]
From Proposition~\ref{prop:7} we get that there are holomorphic
functions $B_1,\dots,B_n$ defined on $\Bbb D$ such that
\begin{multline}
B_m^n(z_m)-B_m^{n-1}(z_m)f_1(z)+\dots+(-1)^{n-1}B_m(z_m)f_{n-1}(z)+
(-1)^nf_n(z)\equiv0\\
z=(z_1,\dots,z_m)\in\DD^n,\ m=1,\dots,n.
\end{multline}
Hence, $f=\pi\big(B_1,\dots,B_n)$. So,
$B=(B_1,\dots,B_n):\DD^n\to\DD^n$ is a proper holomorphic mapping.
From this we have the proof.
\end{proof}

\section{The Shilov boundary}

%

We start with the description of the Shilov boundary of $\GG_n$.
The description of the Shilov boundary in the special case $n=2$
can be found in \cite{Jar-Pfl4} (see also \cite{Agl-You2}).

\begin{lemma} The set $\delta_n$ is the Shilov boundary of
$\GG_n$.
\end{lemma}
\begin{proof} It is easy to see that the modulus of any function from
$\cC(\bar\GG_n)\cap\OO(\GG_n)$ attains its maximum in
$\delta_n$. To finish the proof it is sufficient to show that for
any $z^0\in\delta_n$ there is a function
$F\in\cC(\bar\GG_n)\cap\OO(\GG_n)$ such that $|F|$ attains its
strict maximum at $z^0$.

Fix $z^0:=\pi_n(\lambda_1^0,\ldots,\lambda_n^0)\in\delta_n$,
$|\lambda_1^0|=\ldots=|\lambda_n^0|=1$. Let $B_1$ be a finite
Blaschke product such that
$B_1(\lambda_1^0)=\ldots=B_1(\lambda_n^0)=1$ (see \cite{Abr-Fis},
\cite{You}). Let $B_1^{-1}(1)=\{\lambda_1^0,\ldots,\lambda_m^0\}$,
where $m\geq n$. Let $B_2$ be a finite Blaschke product such that
$B_2(\lambda_1^0)=\ldots=B_2(\lambda_n^0)=1$ and
$B_2(\lambda_{n+1}^0)=\ldots=B_2^0(\lambda_{m})=-1$ (use once more
\cite{Abr-Fis}, \cite{You}). Let $B:=B_1+B_2$. Note that $|B|\leq
2$ on $\bar\DD$ and
$B^{-1}(2)=\{\lambda_1^0,\ldots,\lambda_n^0\}$. Define
$$F(z):=1+\sum\sb{j=1}\sp{n}B(\lambda_j),
$$
where $z=\pi(\lambda_1,\ldots,\lambda_n)\in\bar\GG_n$,
$\lambda_1,\ldots,\lambda_n\in\bar\DD$. Then
$F\in\cC(\bar\GG_n)\cap\OO(\GG_n)$ and $|F|$ attains its strict
maximum (equal to $1+2n$) at $z^0$.
\end{proof}

\section{The Bergman kernel}

It easily follows from the properties of proper holomorphic
mappings that
\begin{multline} \mathcal J_n:=\{\lambda\in\CC^n:\det
\pi^{\prime}(\lambda)=0\}=\\
\{\lambda\in\CC^n:\lambda_j=\lambda_k\text{ for some $j\neq k$,
$j,k=1,\ldots,n$.}\}.
\end{multline}
Denote by $K_D$ the Bergman kernel of the domain $D\subset\CC^n$
(see e.g. \cite{Jar-Pfl3}).
\begin{prop}\label{prop:9}
$$
K_{\GG_n}(\pi_n(\lambda),\pi_{n}(\mu))=
\frac{\det\left[\frac{1}{(1-\lambda_j\bar\mu_k)^2}\right]_{1\leq
j,k\leq
n}}{\pi^n\det\pi_n^{\prime}(\lambda)\det\overline{\pi_n^{\prime}(\mu)}}
$$
for any $\lambda,\mu\in\DD_n\setminus\mathcal J_n$.
\end{prop}

Let $\Sigma_n$ denote the group of all permutations of the set
$I_n$. For $\sigma\in\Sigma_n$,
$\lambda=(\lambda_1,\ldots,\lambda_n)\in\CC^n$ denote
$\lambda_{\sigma}:=(\lambda_{\sigma(1)},\ldots,\lambda_{\sigma(n)})$.

\begin{proof} From the formula for the Bergman kernel of the
polydisc and from the formula for the behavior of the Bergman
kernel under proper holomorphic mappings (see \cite{Bell}) we get
for any $\lambda,\mu\in\DD_n\setminus\mathcal J_n$
\begin{multline} K_{\GG_n}(\pi_n(\lambda),\pi_{n}(\mu))=
\frac{1}{\det\pi_n^{\prime}(\lambda)}\sum\sb{\sigma\in\Sigma_n}
K_{\DD^n}(\lambda,\mu_{\sigma})\frac{1}{\det\overline{\pi_n^{\prime}(\mu_{\sigma})}}
=\\
\frac{1}{\pi^n\det\pi_n^{\prime}(\lambda)\det\overline{\pi_n^{\prime}(\mu)}}
\sum\sb{\sigma\in\Sigma_n}\frac{(-1)^{\sgn\sigma}}{(1-\lambda_j\bar\mu_{\sigma(k)})^2}=
\frac{\det\left[\frac{1}{(1-\lambda_j\bar\mu_k)^2}\right]_{1\leq
j,k\leq
n}}{\pi^n\det\pi_n^{\prime}(\lambda)\det\overline{\pi_n^{\prime}(\mu)}}.
\end{multline}
\end{proof}

The formula above extends analytically to a formula on
$\GG_n\times\GG_n$. It would be interesting to find a more handy
formula for $K_{\GG_n}$. Below we deliver such a formula in the
case $n=2$. We start with the simplification of the denominator in
the formula for $K_{\GG_n}$.

\begin{lemma}\label{lem:10}
$\det\pi_1'(\lambda)=1$, $\lambda\in\CC$, and
$\det\pi_n^{\prime}(\lambda)=\prod_{1\leq j<k\leq
n}(\lambda_j-\lambda_k)$,
$\lambda=(\lambda_1,\dots,\lambda_n)\in\CC^n$, for any $n\ge2$.
\end{lemma}

\begin{proof} We prove by induction on $n$. For $n=2$ we get it by easy computations.
Put $\pi_{n,0}\equiv1$ for any $n\ge1$. Note that
\begin{equation}
\dif{\pi_{n,k}}{\lambda_1}(\lambda_1,\dots,\lambda_n)-
\dif{\pi_{n,k}}{\lambda_n}(\lambda_1,\dots,\lambda_n)=
(\lambda_n-\lambda_1)
\dif{\pi_{n-1,k-1}}{\lambda_1}(\lambda_1,\dots,\lambda_{n-1}).
\end{equation}
Similar equation holds for any pair $(\lambda_j,\lambda_n)$. From
this we get
\begin{multline}
\det\pi_n^{\prime}(\lambda_1,\dots,\lambda_n)= \det\Big[
\dif{\pi_{n,k}}{\lambda_j}(\lambda_1,\dots,\lambda_n)\Big]_{j,k=1,\dots,n}\\
=(\lambda_1-\lambda_n)\dots(\lambda_{n-1}-\lambda_n) \det\Big[
\dif{\pi_{n-1,k}}{\lambda_j}(\lambda_1,\dots,\lambda_{n-1})\Big]_{j,k=1,\dots,n-1}\\
= (\lambda_1-\lambda_n)\dots(\lambda_{n-1}-\lambda_n)
\det\pi_{n-1}^{\prime}(\lambda_1,\dots,\lambda_{n-1}).
\end{multline}
\end{proof}

%
%
%

In the case $n=2$ elementary calculation shows.
\begin{prop}
\begin{multline}\label{eq:25}
K_{\GG_2}(\pi_2(\lambda),\pi_2(\mu))=\frac{2-(\bar\mu_1+\bar\mu_2)(\lambda_1+\lambda_2)+
2\lambda_1\lambda_2\bar\mu_1\bar\mu_2}{\pi^2((1-\lambda_1\bar\mu_1)(1-\lambda_2\bar\mu_2)(1-\lambda_1\bar\mu_2)
(1-\lambda_2\bar\mu_2))^2},\\
\lambda,\mu\in\DD_2.
\end{multline}
In particular, $\GG_2$ is the Lu-Qi-Keng domain, i.e.
$K_{\GG_2}(z,w)\neq 0$ for any $z,w\in\GG_2$.
\end{prop}

\begin{proof} To get the desired formula it is sufficient to apply
the formula from Proposition~\ref{prop:9}, Lemma~\ref{lem:10}, and
then make elementary calculations. To prove that the domain
$\GG_2$ is Lu Qi-keng, it is sufficient, in view of the form of
the automorphisms of $\GG_2$, to verify that
$K_{\GG_2}(\pi_2(\lambda),\pi_2(\mu))\neq 0$,
$\lambda,\mu\in\DD^2$ under the additional assumption $\mu_2=0$.
But this easily follows from the obtained formula \eqref{eq:25}.
\end{proof}

\begin{remark} It would be interesting to find a more
effective formula for the Bergman kernel for $\GG_n$, $n\geq
3$. Moreover, the problem whether the domain $\GG_n$, $n\geq
3$, is Lu Qi-keng, is open, too.
\end{remark}

\section{The spectral ball}

Define $\Omega_n:=\{W\in\mathcal M_n(\CC):r(W)<1\}$, where $r(W)$
denotes the spectral radius in $\mathcal M_n(\CC)$.

Denote also the following mapping
$$
\Psi_n:\Omega_n\owns W\mapsto\pi_n(\sigma(W))\in\GG_n,
$$
where $\sigma(W)$ denotes the spectrum of $W$. Note that $\Psi_n$
is a holomorphic mapping, which is onto but not one-to-one.

\begin{lemma}\label{lem:13} Let $W\in\Omega_n$. Put
$\sigma(W):=\{\lambda_1,\ldots,\lambda_n\}$,
$\lambda_1,\ldots,\lambda_n\in\CC$ (it may happen that
$\lambda_j=\lambda_k$ for some $j\neq k$). Then there is a
holomorphic mapping $f:\CC\mapsto \mathcal M_n(\CC)$ such that
$f(0)=W$, $f(1)=\diag(\lambda_1,\ldots,\lambda_n)$ and
$\sigma(f(\lambda))=\{\lambda_1,\ldots,\lambda_n\}$ for any
$\lambda\in\CC$.
\end{lemma}

\begin{proof} We proceed as in \cite{Cos}. There is an invertible
matrix $X\in\mathcal M_n(\CC)$ and an upper triangular matrix
$$
S=\bmatrix\lambda_1 & x_{1,2} & x_{1,3} & \ldots & x_{1,n-1} & x_{1,n}\\
          0       &\lambda_2 & x_{2,3} & \ldots & x_{2,n-1} & x_{2,n}\\
          \ldots & \ldots & \ldots & \ldots & \ldots & \ldots\\
          0      & 0 & 0 & \ldots &\lambda_{n-1} & x_{n-1,n}\\
          0     & 0 & 0 & 0 & 0 & \lambda_n
          \endbmatrix
          $$
such that $W=XSX^{-1}$. Then we can find a matrix $Y\in\mathcal
M_n(\CC)$ such that $X=e^Y$ (see e.g.~\cite{Au}) and, therefore,
$W=e^YSe^{-Y}$. Consider the holomorphic mappings
$g_{j,k}:\CC\mapsto\CC$, $j=1,\ldots,n-1$, $k=j+1,\ldots,n$ such
that $g_{j,k}(0)=0$ and $g(1)=x_{j,k}$. Then we define
$f:\CC\mapsto\Omega_n$ as follows
$$
f(\lambda):=e^{\lambda Y} \bmatrix\lambda_1 & g_{1,2}(\lambda) &
g_{1,3}(\lambda)
& \ldots & g_{1,n-1}(\lambda) & g_{1,n}(\lambda)\\
          0       &\lambda_2 & g_{2,3}(\lambda) & \ldots & g_{2,n-1}(\lambda) & g_{2,n}(\lambda)\\
          \ldots & \ldots & \ldots & \ldots & \ldots & \ldots \\
          0      & 0 & 0 & \ldots &\lambda_{n-1} & g_{n-1,n}(\lambda)\\
          0     & 0 & 0 & \ldots & 0 & \lambda_n
          \endbmatrix e^{-\lambda Y}.
          $$
\end{proof}

\begin{prop}\label{prop:14} Let $F:\Omega_n\mapsto\Omega_n$ be
holomorphic. Then there is a holomorphic mapping $\tilde F:\GG_n\mapsto\GG_n$ such that
\begin{equation}
\tilde F(\Psi_n(W))=\Psi_n(F(W)),\;W\in\mathcal M_n(\CC).
\end{equation}
\end{prop}

\begin{proof} It is sufficient to show that for any $W\in\Omega_n$
\begin{equation*}
\sigma(F(W))=\sigma(F(\diag(\lambda_1,\ldots,\lambda_n))),
\end{equation*}
where
$\sigma(W)=\{\lambda_1,\ldots,\lambda_n\}$. In view of the
previous lemma there is a holomorphic mapping $f:\CC\mapsto\Omega_n$ such that $f(0)=W$,
$f(1)=\diag(\lambda_1,\ldots,\lambda_n)$.

The mapping $\Psi_n\circ F\circ f:\CC\mapsto\GG_n$ is
holomorphic. Since $\GG_n$ is bounded we get that $\Psi_n\circ
F\circ f$ is constant. In particular,
$\pi_n(\sigma(\diag(\lambda_1,\ldots,\lambda_n)))=\pi_n(\sigma(F(W))$,
which finishes the proof.
\end{proof}

For a finite set $\emptyset\ne P\subset D$ denote $g_D(P,\cdot)$
the pluricomplex Green function with the pole set in $P$ (see
e.g.~\cite{Edi-Zwo}). If $P=\{p\}$ then we put
$g_D(p,\cdot)=g_D(\{p\},\cdot)$.

\begin{prop}\label{prop:15} Let $W_1,W_2\in \Omega_n$. Then
\begin{equation*}
g_{\Omega_n}(W_1,W_2)=g_{\Omega_n}(W_1,\diag(\lambda_1,\ldots.\lambda_n)),
\end{equation*}
where $\sigma(W_2)=\{\lambda_1,\ldots,\lambda_n\}$.
\end{prop}

\begin{proof} Let $f:\CC\mapsto\Omega_n$ be such that
$f(0)=W_2$, $f(1)=\diag(\lambda_1,\ldots,\lambda_n)$ (use
Lemma~\ref{lem:13}). Since the function
$u:=g_{\Omega_n}(W_1,f(\cdot)): \CC\mapsto[-\infty,0)$ is
subharmonic, $u$ must be constant. Therefore,
$g_{\GG_n}(W_1,W_2)=u(0)=u(1)=g_{\Omega_n}(W_1,\diag(\lambda_1,\ldots,\lambda_n))$.
\end{proof}

A domain $D\subset\CC^n$ is called hyperconvex if there exists a
negative plurisubharmonic exhaustion $u$ of $D$, i.e.~$\{z\in D:
u(z)<-\epsilon\}\Subset D$ for any $\epsilon>0$ (see
e.g.~\cite{K}). Note that any hyperconvex domain is pseudoconvex.

We have the following result.
\begin{prop}\label{prop:16} $\GG_n$ is a hyperconvex domain for any $n\ge 1$.
\end{prop}
\begin{proof} We know that $\GG_n=\pi_n(\DD^n)$. Let
$q(z)=\max\{|z_1|,\dots,|z_n|\}$, where
$z=(z_1,\dots,z_n)\in\CC^n$. Put $u(w)=\max\log q(\pi_n^{-1}(w))$,
where $w\in\CC^n$. Note that $u_{|\GG_n}$ is a negative
plurisubharmonic exhaustion of $\GG_n$ (use Proposition 2.9.26 in
\cite{K}).
\end{proof}

The following result is a partial generalization of the main
result in \cite{Ran-W}.
\begin{thm} Let $F:\Omega_n\mapsto\Omega_n$ be a proper
holomorphic mapping. Then there is a finite Blaschke product $B$
such that $\sigma(F(W))=B(\sigma(W))$, $W\in \Omega_n$.
\end{thm}

\begin{proof} Let $\tilde F:\GG_n\mapsto\GG_n$ be as in
Proposition~\ref{prop:14}. Because of Theorem~\ref{thm:1} it is
sufficient to show that $\tilde F$ is proper. It is sufficient to
show that for any sequence $\GG_n\supset(z^{\nu})$ such that
$z^{\nu}\to\partial\GG_n$ there is a subsequence $(z_{\nu_k})$
such that $\tilde F(z^{\nu_k})\to\partial\GG_n$. Fix such a
sequence $(z^{\nu})$. Let $(W^{\nu})\subset\Omega_n$ be such a
sequence that $\Psi_n(W^{\nu})=z^{\nu}$.  Let $Q\in\Omega_n$ be
any point such that $Q\not\in F(\{\det F^{\prime}=0\})$. Let $N$
denote the multiplicity of the proper mapping $F$. Denote
$\{P_1,\ldots,P_N\}=F^{-1}(Q)$. Denote
$\sigma(F(W^{\nu})):=\{\lambda_1^{\nu},\ldots,\lambda_n^{\nu}\}$.
Then because of Proposition~\ref{prop:15} and the behaviour of the
Green function under proper holomorphic mappings (see
\cite{Edi-Zwo}) we get
\begin{multline}
g_{\Omega_n}(Q,\diag(\lambda_1^{\nu},\ldots,\lambda_n^{\nu}))=
g_{\Omega_n}(Q,F(W^{\nu}))=\\
g_{\Omega_n}(\{P_1,\ldots,P_N\}, W^{\nu})\geq
\sum_{j=1}^{N}g_{\Omega_n}(P_j,W^{\nu})\geq\sum_{j=1}^Ng_{\GG_n}(\Psi_n(P_j),z^{\nu}).
\end{multline}
Now the hyperconvexity of $\GG_n$ and the convergence
$z^{\nu}\to\partial\GG_n$ imply that the last expression tends to
$0$ as $\nu\to\infty$ (see e.g.~\cite{K}).

Choosing a subsequence, if necessary, we may assume that
$\lambda_j^{\nu}\to\lambda_j$, $j=1,\ldots,n$, where
$|\lambda_j|\leq 1$. Now the assumption $|\lambda_j|<1$,
$j=1,\ldots,n$, would imply, because of the upper-semicontinuity
of the Green function, that
$g_{\Omega_n}(Q,\diag(\lambda_1,\ldots,\lambda_n))=0$ --
contradiction. Consequently, $\tilde
F(z^{\nu})=\Psi_n(W^{\nu})=\Psi_n(\diag(\lambda_1^{\nu},\ldots,\lambda_n^{\nu}))\to
\Psi_n(\diag(\lambda_1,\ldots,\lambda_n))\in\partial\GG_n$.
\end{proof}

%
%
%
%

\bibliographystyle{amsplain}


\end{document}